\documentclass[11pt]{article}
\usepackage{amsfonts,epsf,amsmath,amssymb}
\usepackage{epsfig}
\usepackage{latexsym}
\usepackage{pgf,tikz}

\newtheorem{theorem}{\bf Theorem}[section]
\newtheorem{corollary}[theorem]{\bf Corollary}
\newtheorem{lemma}[theorem]{\bf Lemma}
\newtheorem{proposition}[theorem]{\bf Proposition}

\newtheorem{question}[theorem]{\bf Question}

\newcommand{\proof}{\noindent{\bf Proof.\ }}
\newcommand{\qed}{\hfill $\square$ \bigskip}

\begin{document}

\title{On disjoint hypercubes in Fibonacci cubes}

\author{
Sylvain Gravier \footnote{CNRS/UJF - Institut Fourier - FR Maths �  Modeler, 100 rue des maths - BP 74, 38402 Saint Martin d'Hères}
\and 
Michel Mollard\footnote{Institut Fourier, CNRS Universit\'e Joseph Fourier, email: michel.mollard@ujf-grenoble.fr}
\and 
Simon \v Spacapan\footnote{University of Maribor, FME, email: simon.spacapan@um.si.}
\footnote{This work is supported by ministry of Education of Slovenia, grants P1-0297 and J1-4106.}
\and 
Sara Sabrina Zemlji\v c\footnote{Institute of Mathematics, Physics and Mechanics, Ljubljana, email: sara.zemljic@gmail.com}
}
\date{\today}
\maketitle

\begin{abstract}
\noindent The {\em Fibonacci cube} of dimension $n$, denoted as $\Gamma_n$,  is the subgraph of  $n$-cube $Q_n$ induced by vertices with no consecutive 1's. We study the maximum number of disjoint subgraphs in $\Gamma_n$ isomorphic to $Q_k$, and denote this number by $q_k(n)$. We prove several recursive results for $q_k(n)$, in particular we prove that 
$q_{k}(n) = q_{k-1}(n-2) + q_{k}(n-3)$. We also prove a closed formula in which $q_k(n)$ is given in terms of  Fibonacci numbers, and finally we give the generating function for the sequence  $\{q_{k}(n)\}_{n=0}^{ \infty}$.

\end{abstract}


\noindent
{\bf AMS Subj. Class. (2010)}: 

\section{Introduction}

Fibonacci cubes have been studied  in several contexts, and the aim of their study goes beyond purely theoretical considerations. 
In mathematical chemistry the concept is related to perfect matchings in hexagonal graphs, and on the other hand,  the structure of perfect matchings in these graphs is used, in theoretical chemistry, to determine the stability of benzenoid molecules. In this context it was proved that the resonance graphs of fibonacenes are isomorphic to Fibonacci cubes \cite{san-pet, kemija1, kemija2}. In computer science  Fibonacci cubes are interesting from algorithmic point of view, since many algorithms that are known to be polynomial (or even linear) on the class 
of hypercubes work  as well on Fibonacci cubes \cite{hsu,cong}. 

Besides many appealing applications of Fibonacci cubes, these graphs have also been studied from a purely theoretical point of view. 
One of the reasons is due to their nice recursive structure and properties derived from it (see \cite{survey}). In particular we mention that several graph theoretical invariants are easily determined using their structural properties. It was proved in \cite{cong} that every Fibonacci cube has a hamiltonian path, and in \cite{36}  the independence number of Fibonacci cubes is determined. 
Various enumerative sequences of these graphs have been determined as well: number of vertices of a given degree \cite{klmope-2011}, number of vertices of a given eccentricity \cite{camo-2012} and number of isometric subgraphs isomorphic to  $Q_k$ 
 \cite{KlavzarCube}. The counting polynomial of this last sequence is known as cubic polynomial and has very nice properties. 
In  \cite{mollard} the author gave the number of maximal hypercubes of dimension $k$ in $\Gamma_n$  (that is,  induced subgraphs $H$  of $\Gamma_n$ isomorphic to  $Q_k$, such that their exists no induced subgraph $H'$  isomorphic to $Q_{k+1}$, with $H\subseteq H'$) and determined the exact number of maximal hypercubes for every $k$ and $n$.

In this paper we study the maximum number of disjoint subgraphs isomorphic to $Q_k$ in $\Gamma_n$. This problem is relevant in the context of parallel or distribute systems. Indeed we can imagine different tasks occurring in parallel on each $Q_k$ without conflict.
It turns out that the problem has a natural solution which is given by several recursive relations and a generating function (in the following section). In the rest of the introduction we give the definitions and notations.

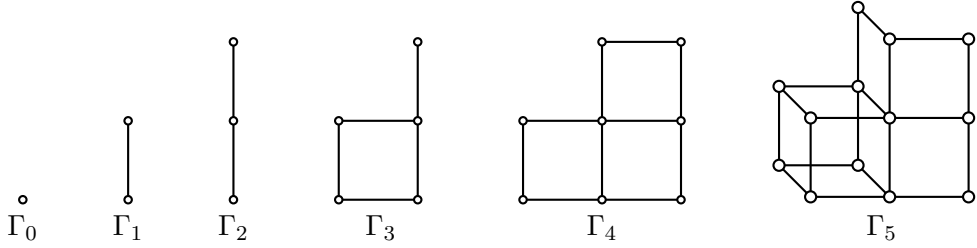
\begin{figure}[ht!] 
		\begin{center}
			\begin{tikzpicture}[scale=0.7,style=thick]
					\def \vr{2pt}
				\draw (0,0) [fill=white] circle (\vr);
					\draw (0,-0.1) [anchor=north] node {$\Gamma_0$};
				\draw (2,0)-- (2,1.5);
				\draw (2,0) [fill=white] circle (\vr);
				\draw (2,1.5) [fill=white] circle (\vr);
					\draw (2,-0.1) [anchor=north] node {$\Gamma_1$};
				\draw (4,0)-- (4,1.5)--(4,3);
				\draw (4,0) [fill=white] circle (\vr);
				\draw (4,1.5) [fill=white] circle (\vr);
				\draw (4,3) [fill=white] circle (\vr);
					\draw (4,-0.1) [anchor=north] node {$\Gamma_2$};
				\draw (7.5,1.5)--(6,1.5)--(6,0)--(7.5,0)-- (7.5,1.5)--(7.5,3);
				\draw (6,0) [fill=white] circle (\vr);
				\draw (6,1.5) [fill=white] circle (\vr);
				\draw (7.5,0) [fill=white] circle (\vr);
				\draw (7.5,1.5) [fill=white] circle (\vr);
				\draw (7.5,3) [fill=white] circle (\vr);
					\draw (6.8,-0.1) [anchor=north] node {$\Gamma_3$};
				\draw (11,0)--(12.5,0)--(12.5,1.5)--(11,1.5)--(11,3)--(12.5,3)--(12.5,1.5);
				\draw (11,1.5)--(11,0)--(9.5,0)--(9.5,1.5)--cycle;
				\draw (9.5,0) [fill=white] circle (\vr);
				\draw (9.5,1.5) [fill=white] circle (\vr);
				\draw (11,0) [fill=white] circle (\vr);
				\draw (11,1.5) [fill=white] circle (\vr);
				\draw (11,3) [fill=white] circle (\vr);
				\draw (12.5,0) [fill=white] circle (\vr);
				\draw (12.5,1.5) [fill=white] circle (\vr);
				\draw (12.5,3) [fill=white] circle (\vr);
					\draw (11,-0.1) [anchor=north] node {$\Gamma_4$};
			\end{tikzpicture}
			$\quad$ $\quad$
			\begin{tikzpicture}[scale=1.05,style=thick]
				\def \vr{2pt}
				\draw (2,0)--(2,1)--(1,1)--(1,0)--(2,0)--(3,0)--(3,1)--(2,1)--(2,2)--(3,2)--(3,1);
				\draw (1.6,2.4)--(1.6,1.4)--(1.6,0.4)--(0.6,0.4)--(0.6,1.4)--(1.6,1.4);
				\draw (1,0)--(0.6,0.4);
				\draw (1,1)--(0.6,1.4);
				\draw (2,0)--(1.6,0.4);
				\draw (2,1)--(1.6,1.4);
				\draw (2,2)--(1.6,2.4);
				\draw (1,0) [fill=white] circle (\vr);
				\draw (2,0) [fill=white] circle (\vr);
				\draw (3,0) [fill=white] circle (\vr);
				\draw (1,1) [fill=white] circle (\vr);
				\draw (2,1) [fill=white] circle (\vr);
				\draw (3,1) [fill=white] circle (\vr);
				\draw (2,2) [fill=white] circle (\vr);
				\draw (3,2) [fill=white] circle (\vr);
				\draw (0.6,0.4) [fill=white] circle (\vr);
				\draw (1.6,0.4) [fill=white] circle (\vr);
				\draw (0.6,1.4) [fill=white] circle (\vr);
				\draw (1.6,1.4) [fill=white] circle (\vr);
				\draw (1.6,2.4) [fill=white] circle (\vr);
				\draw (1.9,-0.1) [anchor=north] node {$\Gamma_5$};
			\end{tikzpicture}
			\caption{Fibonacci cubes $\Gamma_0$, $\Gamma_1$, $\Gamma_2$, $\Gamma_3$, 
				$\Gamma_4$, and $\Gamma_5$}
			\label{fig:cubes}
		\end{center}
	\end{figure}

Let $n$ be a positive integer and denote $[n]=\{1,\ldots, n\}$, and $[n]_0=\{0,\ldots, n-1\}$. 
The {\em $n$-cube}, denoted by $Q_n$, is the graph with vertex set$$V(Q_n)=\{x_1x_2\ldots x_n \, | \, x_i\in [2]_0 ~\text{for}~ i\in [n]\}\,,$$ where 
two vertices 
are adjacent in $Q_n$ if the corresponding strings differ in exactly one position. The {\em Fibonacci $n$-cube}, denoted by $\Gamma_n$, is the subgraph of $Q_n$ induced by vertices with no consecutive 1's. Moreover we define $Q_0 = \Gamma_0$ as a one-vertex graph (see Fig.~\ref{fig:cubes}). 

The vertices of $\Gamma_n$ beginning with 0 induce a subgraph isomorphic to $\Gamma_{n-1}$, and similarly, the vertices beginning with 10 induce a subgraph isomorphic to  $\Gamma_{n-2}$. Denote these subgraphs with $0\Gamma_{n-1}$ and  $10\Gamma_{n-2}$, respectively. We  generalize this notation as follows: a subgraph of $\Gamma_n$, induced by vertices with a common prefix $b0$, $b\in[2]_0^{i-1}$, $i\in[n]$, is denoted by $b0\Gamma_{n-i}$.
Note that there is a matching from $0\Gamma_{n-1}$  to $10\Gamma_{n-2}$,  where each vertex $00x$ of $0\Gamma_{n-1}$ is paired with $10x$ of  $10\Gamma_{n-2}$. We will denote the subgraph of $\Gamma_n$ induced by endvertices of this matching with $00\Gamma_{n-2} \equiv 10\Gamma_{n-2}$. So we may think that $\Gamma_n$ is constructed from $0\Gamma_{n-1}$ and $10\Gamma_{n-2}$, where the only edges between these two subcubes induce a matching that saturates $10\Gamma_{n-2}$. 
Finally we  mention that the recursive structure of Fibonacci cubes gives the following equations 
$$  |\Gamma_0| = 1\,,\ |\Gamma_1| = 2 \ \text{ and } \ \ |\Gamma_n| = |\Gamma_{n-1}| + |\Gamma_{n-2}|,$$
which in turn gives us the Fibonacci numbers $\{F_{n}\}_{n\ge 0}$, i.\ e.\ $|\Gamma_n| = F_{n+2}$.

In the next section we study and determine the maximum number of disjoint subgraphs isomorphic to $Q_k$ of the Fibonacci cube $\Gamma_n$. Denote this number by $q_k(n)$.

\section{Recursive relations}
	\indent 
	In the next table, we summarize the first values of $q_{k}(n)$ obtained by direct inspection of $\Gamma_{n}$.
		\begin{table}[htb]
			\begin{center}
				\begin{tabular}{|c|c c c c c c|}
			     	\hline
					\textit{n} & 0 & 1& 2& 3 &  4 & 5 \\\hline
					$|\Gamma_n|$ & 1 & 2 & 3 & 5 & 8 & 13 \\
					$q_{1}(n)$ & 0 & 1 & 1 & 2 & 4 & 6 \\
					$q_{2}(n)$ & 0 & 0 & 0 & 1 & 1 & 2 \\
					$q_{3}(n)$ & 0 & 0 & 0 & 0 & 0 & 1 \\\hline
				\end{tabular}
				\caption{First values of $q_k(n)$ for $n\in[6]_0$ and $k\in[4]_0$}
				\label{tab:initial}
			\end{center}
		\end{table}

	In the sequel we will  differentiate two cases for $n$, depending on the residue when dividing by 3. Therefore let $r(n)\in [3]_0$ be such that $n\equiv r(n) \,(\mathrm{mod}\, 3)$. By the recursive structure of Fibonacci cubes described in the introduction, the cube $\Gamma_n$ can be partitioned into subgraphs $00\Gamma_{n-2} \equiv 10\Gamma_{n-2}$ and $010\Gamma_{n-3}$. Note that there are also some edges between these two subgraphs, namely the matching between $000\Gamma_{n-3}$ and $010\Gamma_{n-3}$, but in the construction below we shall not need them. Denote this partition by $\dot{\cup}$. In particular
	\begin{equation}
		\label{eq:gamma_recc}
		\Gamma_n = (00\Gamma_{n-2}\equiv 10\Gamma_{n-2})\ \dot{\cup}\ 010\Gamma_{n-3} \,.
	\end{equation}
		$$$$
	Repeating this gives us a complete partitioning of the Fibonacci cube: 
		\begin{align}
			\label{eq:gamma_total_recc}
			\Gamma_n= \dot{\bigcup}_{i=1}^{\lfloor \frac{n}{3}\rfloor} \big( (010)^{i-1} 00\Gamma_{n+1-3i}&\equiv (010)^{i-1} 10\Gamma_{n+1-3i} \big) \dot{\cup}  (010)^{\frac{n-r(n)}{3}} \Gamma_{r(n)}\,.
		\end{align}

	First let us derive the exact size of a maximum matching in a Fibonacci cube. 

	\begin{lemma}
		\label{lem-q_1,n}
		For every $n\in\mathbb{N}$ we have
			\begin{equation*}
						q_{1}(n)=\frac{\vert \Gamma_n\vert - \gamma(n)}{2}
							=\left\lfloor \frac{F_{n+2}}{2} \right\rfloor\,.
				\end{equation*}
		Where $\gamma(n)=0$ if $r(n)=1$, and $\gamma(n)=1$ otherwise. 
	\end{lemma}
		
	\proof
		By~\eqref{eq:gamma_total_recc} each 
			$\big( (010)^{i-1} 00\Gamma_{n+1-3i}\equiv (010)^{i-1} 10\Gamma_{n+1-3i} \big)$ 
		has a perfect matching, and therefore we only need to consider $r(n)$.
		If $r(n)$ is $0$ or $2$ then $\Gamma_{r(n)}$ has a maximum matching that misses one vertex. Otherwise $\Gamma_{r(n)}=\Gamma_{1} $ has a perfect matching. From this we get the lower bound. The upper bound follows from the parity of the number of vertices of Fibonacci cubes: $|\Gamma_n|=F_{n+2}$ it is even if $r(n)=1$ and odd otherwise.
	\qed

	Note that in the case $r(n)=1$ the Fibonacci cube $\Gamma_n$ has a perfect matching. In all the other cases there is only one vertex not covered by a maximum matching. By similar construction we are able to prove even more general case, i.\ e.\ when $k\geq 2$.

	\begin{theorem}
		\label{thm-q_k,n}
		For every $n\geq 3$ and $k\geq 2$
			\begin{equation}
				\label{eq-q_n,k}
				q_{k}(n) = q_{k-1}(n-2) + q_{k}(n-3)\,.
			\end{equation}
	\end{theorem}

	\proof
		Note that \eqref{eq-q_n,k} is equivalent to 
			$$ q_{k}(n) = \sum_{i=1}^{\lfloor \frac{n}{3}\rfloor} q_{k-1}(n+1-3i)\,. $$
		This is indeed true, since for $n\in [3]_0$, $q_{k}(n) = 0$ for any $k\ge 2$ (see Table~\ref{tab:initial}).

		The lower bound follows directly from \eqref{eq:gamma_recc}. We construct a family of disjoint $k$-hypercubes by taking the maximum number of disjoint $(k-1)$-cubes in $00\Gamma_{n-2}$ and connecting each of them with the corresponding $(k-1)$-cube in $10\Gamma_{n-2}$. By doing so, we get $q_{k-1}(n-2)$ $k$-cubes in $00\Gamma_{n-2}\equiv10\Gamma_{n-2}$. Then we add the maximum number of $k$-cubes in $010\Gamma_{n-3}\cong\Gamma_{n-3}$, this is additional $q_{k}(n-3)$ $k$-cubes.

		To prove the upper bound consider the partition of $\Gamma_n$ described in \eqref{eq:gamma_total_recc}.
		Let $\mathcal{U}$ be a family of disjoint $k$-cubes in $\Gamma_n$. We claim that 
			\begin{equation}
				\label{eq:claim1}
				|\cup \mathcal{U} \, \cap \, V\big( (010)^{i-1} 00\Gamma_{n+1-3i}\equiv (010)^{i-1} 10\Gamma_{n+1-3i} \big)| \leq 2^{k} \cdot q_{k-1}(n+1-3i)\,.
			\end{equation}
		To prove this note that for every $k$-cube $Q_k \in\mathcal{U}$, we have 
			$$Q_k \cap \big( (010)^{i-1} 00\Gamma_{n+1-3i}\equiv (010)^{i-1} 10\Gamma_{n+1-3i} \big) \in \{\emptyset, Q_{k-1},Q_k\}\,,$$
		for every $i\leq \lfloor \frac{n}{3} \rfloor$.
		Therefore, every $k$-cube in $\mathcal{U}$ intersects $(010)^{i-1} 00\Gamma_{n+1-3i}$ in some number of disjoint $(k-1)$-cubes (every $k$-cube contained in $(010)^{i-1} 00\Gamma_{n+1-3i}$ is regarded as a disjoint union of two $(k-1)$-cubes). This number is at most $q_{k-1}(n+1-3i)$ and hence \eqref{eq:claim1} follows. Thus 
		$$|\cup \mathcal{U} \, \cap \, V(\Gamma_{n})| \leq 2^{k} \cdot \sum_{i=1}^{\lfloor \frac{n}{3} \rfloor} q_{k-1}(n+1-3i)\,,$$
		which completes the proof.
	\qed

\noindent For every $n\geq 2$ and $k\geq 1$ let us define $\delta_k(n)$ as follows: 
$$\delta_k(n)= \binom{\frac{n+k-2}{3}}{k-1}\,,$$  if $r(n+k)=2$, and $\delta_k(n)= 0$ otherwise. We shall use this notation 
throughout the rest of the paper.

	\begin{corollary}
		\label{cor-q_k,n}
		For every $n\geq 2$ and $k\geq 1$ we have 
			\begin{equation}
				\label{eq-q_k,n1}
				q_{k}(n) = q_{k}(n-1) + q_{k}(n-2) +\delta_k(n).
			\end{equation}
	\end{corollary}

	\proof
		The proof is by induction, where the basis are the cases $k=1$ or $n\leq4$:
			\begin{itemize}
				\item[(i)] The case $k=1$ follows directly from Lemma \ref{lem-q_1,n}.
				\item[(ii)] If $n\leq 4$ and $k\geq 3$ we are done since $q_{k}(n)=q_{k}(n-1)=q_{k}(n-2)=0$ and  $\frac{n+k-2}{3}<k-1$.
				\item[(iii)] If $(k,n)= (2,2)$ then  $q_{2}(2)=q_{2}(1)=q_{2}(0)=0$ and we are done.
				\item[(iv)] If $(k,n)= (2,3)$ then equation \eqref{eq-q_k,n1}  is true since  $q_{2}(3)=1, q_{2}(2)=q_{2}(1)=0$ and $\binom{1}{1}=1$.
				\item[(v)] If $(k,n)= (2,4)$ then equation \eqref{eq-q_k,n1}  is true since  $q_{2}(4)=q_{2}(3)=1, q_{2}(2)=0$.
			\end{itemize}
		Let $k\geq2$ and $n\geq5$. Assume that equation \eqref{eq-q_k,n1} is true for all pairs $(k',n') \neq (k,n)$, such that $1\leq k'\leq k$ and $2\leq n'\leq n$ and let us prove the case $(k,n)$.
		We first prove \eqref{eq-q_k,n1} when $n+k \equiv 0$ or $1 \,(\mathrm{mod}\, 3)$. Since 
			\begin{align*}
				& k-1+n-2 \equiv k+n \equiv 0 \textrm{ or } 1 \,(\mathrm{mod}\, 3) \,,\\
				& k+n-3 \equiv k+n \equiv 0\textrm{ or }1 \,(\mathrm{mod}\, 3)  \,, \\
				& \text{for } \ k-1\geq 1 ,\  n-3\geq2
			\end{align*}
		it follows from the hypothesis, that
			\begin{align*}
				& q_{k-1}(n-2)= q_{k-1}(n-3) + q_{k-1}(n-4)\,,\\
				& q_{k}(n-3) = q_{k}(n-4) + q_{k}(n-5)\,.
			\end{align*}
		Then we have 
			\begin{align*}
				q_{k}(n) &= q_{k-1}(n-2) + q_{k}(n-3)\\
					&= q_{k-1}(n-3) + q_{k-1}(n-4) + q_{k}(n-4) + q_{k}(n-5)\\
					&=q_{k}(n-1)+q_{k}(n-2)\,,
			\end{align*}
		where the first and the last equality follow from \eqref{eq-q_n,k}.
		Next we prove \eqref{eq-q_k,n1} when $n+k\equiv 2 \,(\mathrm{mod}\, 3)$. We have 
			\begin{align*}
				q_{k}(n) &= q_{k-1}(n-2) + q_{k}(n-3)\\
					&= q_{k-1}(n-3) + q_{k-1}(n-4) + \binom{\frac{n+k-5}{3}}{k-2} +q_{k}(n-4) + q_{k}(n-5)+\binom{\frac{n+k-5}{3}}{k-1}\\
					&=q_{k}(n-1)+q_{k}(n-2)+\binom{\frac{n+k-2}{3}}{k-1}\,.
			\end{align*}
		So the assertion is proved.
	\qed

	We are now able to prove a closed formula for the number $q_k(n)$.

	\begin{corollary}
		\label{cor-closefq_k,n}
		For every $n\geq 0$ and $k\ge1$
			\begin{equation}
				\label{eq-closefq_k,n}
				q_{k}(n) = \sum_{i=k-1}^{\left\lfloor\frac{n+k-2}{3}\right\rfloor}
						\binom{i}{k-1}F_{n+k-3i-1}.
			\end{equation}
	\end{corollary}
	
	\proof
		For every $k\geq 1$ let
			$$f_k(x)=\sum_{n=0}^{\infty}q_{k}(n)x^n$$
		be the generating function of the sequence $\{q_{k}(n)\}_{n=0}^{ \infty}$.
		 By Corollary \ref{cor-q_k,n} we have
			\begin{equation}
				\label{eq-q_kcom}
				q_{k}(n) = q_{k}(n-1) + q_{k}(n-2) + \delta_k(n)\,,
			\end{equation}
		for every $k\geq 1$ and $n\geq 2$.
		Fixing a $k\geq 1$ and using (\ref{eq-q_kcom}) we obtain
			$$\sum_{n=2}^{\infty}q_{k}(n)x^n =\sum_{n=2}^{\infty}q_{k}(n-1)x^{n} + \sum_{n=2}^{\infty}q_{k}(n-2)x^{n} + \sum_{n=2}^{\infty}\delta_k(n)x^n\,.$$
		Assume first $k\geq2$. Since $q_{k}(0)=q_{k}(1)=0$ and $\delta_k(0)=\delta_k(1)=0$ we have
		$$\sum_{n=0}^{\infty}q_{k}(n)x^{n}=\sum_{n=2}^{\infty}q_{k}(n)x^{n} =x \sum_{n=0}^{\infty}q_{k}(n)x^{n}+ x^2 \sum_{n=0}^{\infty}q_{k}(n)x^{n} + \sum_{n=0}^{\infty}\delta_k(n)x^n.$$
		Therefore
			$$f_k(x)(1-x-x^2)=\sum_{n=0}^{\infty}\delta_k(n)x^n\,,$$
		and so
			\begin{equation}
				\label{eq-q_kcom,5}
				f_k(x)=\left(\frac{1}{1-x-x^2}\right)\left(\sum_{n=0}^{\infty}\delta_k(n)x^n\right)\,.
			\end{equation}
		For $k=1$ we have $q_{1}(1)=1$, $q_{1}(0)=0$, $\delta_1(1)=1$ and $\delta_1(0)=0$. Therefore
			\begin{align*}
				f_1(x)&=\sum_{n=0}^{\infty}q_{1}(n)x^n=q_{1}(0)+q_{1}(1)x+\sum_{n=2}^{\infty}q_{1}(n)x^n\\
						&=x+ \sum_{n=2}^{\infty}q_{1}(n-1)x^{n} + \sum_{n=2}^{\infty}q_{1}(n-2)x^{n} + \sum_{n=2}^{\infty}\delta_1(n)x^n\\
						&=x+ x\sum_{n=0}^{\infty}q_{1}(n)x^{n} + x^2\sum_{n=0}^{\infty}q_{1}(n)x^{n} + \sum_{n=0}^{\infty}\delta_1(n)x^n-\delta_1(1)x\,.
			\end{align*}
		It follows that \eqref{eq-q_kcom,5} is true also for $k=1$. The generating function of the Fibonacci sequence is (see for example \cite{survey}) 
			$$\sum_{n=0}^{\infty}F_nx^n=\frac{x}{1-x-x^2}\,,$$
		and since $F_0=0$, we have $$\frac{1}{1-x-x^2}=\sum_{n=0}^{\infty}F_{n+1}x^n\,.$$
		So by \eqref{eq-q_kcom,5} we have
			$$f_k(x)=\left(\sum_{n=0}^{\infty}F_{n+1}x^n\right ) \cdot \left( \sum_{n=0}^{\infty}\delta_k(n)x^n\right )\,, $$
		and therefore
			\begin{equation}
				\label{eq-q_k,final}
				q_{k}(n) =\sum_{m=0}^{n}{ \delta_k(m)F_{n-m+1}}\,.
			\end{equation}
		Note that for $m \in[n+1]_0$ we have $\delta_k(m)\neq 0$ if and only if
		$ m+k-2 \equiv 0 \,(\mathrm{mod}\, 3)$ and $\frac{m+k-2}{3}\geq k-1$.
		By setting $m+k-2=3i$ we get $3i-k+2\in [n+1]_0$ and $i\geq k-1$.
		So for $m= 3i-k+2$ and $i=k-1,\ldots,\lfloor \frac{n+k-2}{3}\rfloor$ it holds $\delta_k(m)\neq 0$ . Replacing $ 3i-k+2$ for $m$ in equation (\ref{eq-q_k,final}) we obtain the equation (\ref{eq-closefq_k,n}).
	\qed

	By computing the generating function $f_k(x)$ we are able to get additional information on the sequence $\{q_{k}(n)\}_{n=0}^{ \infty}$.

	\begin{corollary}
		\label{cor-gf_k,n}
		For every $k\geq1$ the generating function of the sequence $\{q_{k}(n)\}_{n=0}^{ \infty}$ is
			\begin{equation}
				f_k(x)=\frac{x^{2k-1}}{(1-x^3)^k(1-x-x^2)}\,.
			\end{equation}
	\end{corollary}

	\proof
		First let us prove that 
			$$\frac{1}{(1-x^3)^k}= \sum_{i=0}^{\infty}{ \binom{i}{k-1}x^{3(i-k+1)}}\,.$$
		We proceed by induction on $k$. Base of the induction is clear, since we have 
			$$\frac{1}{1-x^3}= \sum_{i=0}^{\infty}{x^{3i}}\,.$$
		Let now $k\geq2$. Then 
			\begin{align*}
				\frac{1}{(1-x^3)^{k}} &= \left(\sum_{i=0}^{\infty}{ \binom{i}{k-2}x^{3(i-k+2)}}\right)\left(\sum_{j=0}^{\infty}x^{3j}\right)\\
					&=\sum_{i=0}^{\infty}{ \binom{i+1}{k-1}x^{3(i-k+2)}}\\
					&=\sum_{i=0}^{\infty}{ \binom{i}{k-1}x^{3(i-k+1)}}\,,
			\end{align*}
		since summing $\displaystyle \sum_{l=0}^{i}{\binom{l}{k-2}} =\binom{i+1}{k-1}$. 
		Recall that 
			$$\delta_k(n)= \binom{\frac{n+k-2}{3}}{k-1}\,,$$ 
 if $r(n+k)=2$, and $\delta_k(n)= 0$ otherwise. Therefore
		 the generating function of the sequence $\{\delta_k(n)\}_{n=0}^{ \infty}$ equals 
			$$\sum_{i=k-1}^{\infty}{ \binom{i}{k-1}x^{3i-k+2}}=x^{2k-1}\left(\sum_{i=0}^{\infty}{ \binom{i}{k-1}x^{3(i-k+1)}}\right)=\frac{x^{2k-1}}{(1-x^3)^k}\,.$$

		Since by \eqref{eq-q_kcom,5} the sequence $\{q_{k}(n)\}_{n=0}^{ \infty}$ is the convolution of the sequences $\{F_{n+1}\}_{n=0}^{ \infty}$ and $\{\delta_k(n)\}_{n=0}^{ \infty}$ the result follows.

\qed


\section{Concluding remark and open problem}
\begin{proposition}
\label{bound}
		For every  $k\ge1$ and $n\geq 2k-2$
			\begin{equation}
				q_{k}(n) \geq F_{n-2k+2}.
			\end{equation}
\end{proposition}
\proof
Let  $k\ge1$. By corollary 3.3 in \cite{KlavzarCube}, or by our corollary  \ref{cor-closefq_k,n}, we have $	q_{k}(2k-2)=0=F_0$ and $	q_{k}(2k-1)=1=F_1$.
Now let $n\geq 2k$ and assume $	q_{k}(n-1)\geq F_{n-2k+1}$ and $	q_{k}(n-2)\geq F_{n-2k}$. Using corollary  \ref{cor-q_k,n} we obtain
$$	q_{k}(n)\geq q_{k}(n-1)+q_{k}(n-2)\geq F_{n-2k+1}+ F_{n-2k}=F_{n-2k+2}.$$
The result follows by induction.
\qed
\begin{question}
Determine $$\lim_{n\to \infty} \frac{	q_{k}(n) }{|V(\Gamma_n)|}.$$
\end{question}
For any $k\geq 1$,
 $$\frac{q_{k}(n) }{|V(\Gamma_n)|}\geq \frac{	F_{n-2k+2}}{F_{n+2}}\,,$$
by Proposition \ref{bound}.
Therefore $$\lim_{n\to \infty} \frac{	q_{k}(n) }{|V(\Gamma_n)|}\geq \frac{1}{\phi^{2k}}\,,$$ 
where $\phi=(1+\sqrt{5})/2$. This bound is weak since 
$$\lim_{n\to \infty} \frac{	q_{1}(n) }{|V(\Gamma_n)|}=\frac{1}{2}\,.$$
For $k>1$ it will be interesting to study, asymptotically, the ratio of vertices of $\Gamma_n $ covered by a maximum set of disjoint subgraphs isomorphic to $Q_k$, thus compare $\lim_{n\to \infty} \frac{	q_{k}(n) }{|V(\Gamma_n)|}$ and $\frac{1}{2^k}$.


\begin{thebibliography}{99}

\bibitem{camo-2012}
	A.~Castro, M.~Mollard,
	The eccentricity sequences of Fibonacci and Lucas cubes,
	Discrete Math. 312 (2012) 1025--1037.

\bibitem {cong}  
	B.~Cong, S.~Zheng, S.~Sharma, 
	On simulations of linear arrays, rings and 2d meshes on Fibonacci cube networks,
	In: Proceedings of the 7th International Parallel Processing Symphosium, 1993, 747--751.

\bibitem{hsu} 
	W.-J.~Hsu, 
	Fibonacci cubes---a new interconnection technology,
	IEEE Trans. Parallel Distrib.\ Syst. 4 (1993) 3--12.

\bibitem{survey} 
	S.~Klav\v zar, 
	Structure of Fibonacci cubes: a survey, 
	J.\ Comb.\ Optim. 25 (2013) 505--522.

  \bibitem{KlavzarCube}
	S.~Klav\v{z}ar, M.~Mollard,
	Cube polynomial of Fibonacci and Lucas cube,
	Acta Appl.\ Math. 117 (2012) 93--105.

\bibitem{klmope-2011}
	S.~Klav\v{z}ar, M.~Mollard, M.~Petkov\v sek,
	The degree sequence of Fibonacci and Lucas cubes,
 	Discrete Math. 311 (2011) 1310--1322.

\bibitem{san-pet} 
	S.~Klav\v zar, P.~\v Zigert, 
	Fibonacci cubes are the resonance graphs of Fibonaccenes,
	Fibonacci Quart. 43 (2005) 269--276.

\bibitem{mollard}
	M.~Mollard,
	Maximal hypercubes in Fibonacci and Lucas cubes,
	Discrete Appl.\ Math. 160 (2012) 2479--2483.

\bibitem{36} 
	E.~Munarini, N.~Zagaglia Salvi, 
	Structural and enumerative properties of the Fibonacci cubes,
	Discrete Math. 255( (2002) 317--324.

\bibitem{kemija1}
	H.~Zhang, P.~C.~B.~Lam, W.~C.~Shiu, 
	Resonance graphs and a binary coding for the 1-factors of benzenoid systems,
	SIAM J.\ Discrete Math. 22 (2008) 971--984.

\bibitem{kemija2}
	H.~Zhang, L.~Ou, H.~Yao, 
	Fibonacci-like cubes as Z-transformation graphs,
	Discrete Math., 309 (2009) 1284--1293.


  
 
  




\end{thebibliography}
\end{document}